\documentclass[11pt,a4paper,twoside]{article}

\usepackage{amsmath,amsthm,amstext,amscd,amssymb,euscript,mathrsfs,color}
\usepackage{calrsfs}
\usepackage{epsfig}

\newcommand{\R}{\mathbb R}
\newcommand{\N}{\mathbb N}

\newcommand{\E}{\mathbb E}

\newcommand{\epsi}{\ensuremath{\epsilon}}

\newcommand{\loc}{\mathcal{L}}

\def\1{{\mathchoice {\rm 1\mskip-4mu l} {\rm 1\mskip-4mu l}
{\rm 1\mskip-4.5mu l} {\rm 1\mskip-5mu l}}}

\newtheorem{theorem}{{\small T}{\scriptsize HEOREM}}[section]
\newtheorem{corollary}{{\bf{\small C}{\scriptsize OROLLARY}}}[section]
\newtheorem{proposition}{{\bf{\small P}{\scriptsize ROPOSITION}}}[section]
\newtheorem{lemma}{{\bf{\small L}{\scriptsize EMMA}}}[section]
\newtheorem{remark}{{\bf{\small R}{\scriptsize EMARK}}}[section]
\newtheorem{definition}{{\bf{\small D}{\scriptsize EFINITION}}}[section]

\renewenvironment{proof}[1]
{\noindent{{\bf{\small{ P}{\scriptsize ROOF}}}.}\hspace{0.1cm} #1} {$\;\qed$\newline}

\newcommand{\beq}{\begin{eqnarray}}
\newcommand{\eeq}{\end{eqnarray}}

\newcommand{\ba}{\begin{align*}}
\newcommand{\ea}{\end{align*}}

\newcommand{\be}{\begin{equation}}
\newcommand{\ee}{\end{equation}}

\newcommand{\bl}{\begin{lemma}}
\newcommand{\el}{\end{lemma}}

\newcommand{\br}{\begin{remark}}
\newcommand{\er}{\end{remark}}

\newcommand{\bt}{\begin{theorem}}
\newcommand{\et}{\end{theorem}}

\newcommand{\bd}{\begin{definition}}
\newcommand{\ed}{\end{definition}}

\newcommand{\bp}{\begin{proposition}}
\newcommand{\ep}{\end{proposition}}

\newcommand{\bc}{\begin{corollary}}
\newcommand{\ec}{\end{corollary}}

\newcommand{\bpr}{\begin{proof}}
\newcommand{\epr}{\end{proof}}

\newcommand{\bi}{\begin{itemize}}
\newcommand{\ei}{\end{itemize}}

\newcommand{\ben}{\begin{enumerate}}
\newcommand{\een}{\end{enumerate}}

\newcommand{\caE}{{\mathrsfs E}}

\newcommand{\caX}{{\mathcal X}}

\begin{document}
\title{Duality and stationary distributions of wealth distribution models}
\author{
Pasquale Cirillo$^{\textup{{\tiny(a)}}}$,\\
Frank Redig$^{\textup{{\tiny(a)}}}$,\\
Wioletta Ruszel$^{\textup{{\tiny(a)}}}$,
\\
{\small $^{\textup{(a)}}$
Delft University of Technology}, {\small Mekelweg 4 2628 CD Delft , The Netherlands}
\\
}
\maketitle

\pagenumbering{arabic}

\begin{abstract}
We analyze a class of energy and wealth redistribution models.
We characterize their stationary measures and show that they have a discrete
dual process. In particular we show that the wealth distribution model
with non-zero saving propensity can never have invariant product measures.
We also introduce diffusion processes associated to the
wealth distribution models by ``instantaneous thermalization''.
\end{abstract}

\section{Introduction}

Wealth distribution models represent a flourishing field of the econophysical literature, see e.g. \cite{CCCC} for a recent overview. 
They represent simplified models of an economy, where at random instances agents exchange wealth, whereas
the total wealth is conserved. These models are inspired from
kinetic theory, since the exchange of wealth among economic agents is reminiscent of
the exchange of energy in the collisions among the particles in a gas. Therefore, these models
have striking similarities with the stochastic processes of energy exchange that are used
as microscopic toy models of ``the Fourier law'' of heat conduction.\\
One particular aspect of wealth distribution models which makes them quite different
from analogous models in mathematical physics, such as the KMP model red for Fourier law \color{black} \cite{KMP}
or the energy and particle transport models studied in \cite{GKRV}, is the presence of  the so-called ``propensity to save'', i.e.
the tendency of agents to have only a fraction of their wealth involved in transaction events, while the other amount constitutes the savings, which are not transacted and allow for the accumulation of wealth.\\
It is therefore interesting to study wealth distribution models from the point of view of mathematical
physics, using in particular the techniques of duality developed in \cite{GKRV} and \cite{CGGR}. The aim of this paper is to present first results in that direction, by restricting our attention to a general class of wealth distribution models on two sites. \\

For this class of models, we characterize their invariant measures and we prove that for non-zero propensity,
there does not exist an invariant product measures, i.e., non-zero propensity necessarily leads to dependencies
between the wealth of agents. We also show that these models always have a discrete dual process, but the duality function factorizes if and only if the propensity to save is zero and the redistribution measure a beta distribution, corresponding to the unique
case where we do have product invariant measures. In other words, factorization can only take place when the wealth distribution model is nothing but an energy redistribution model of the type described in \cite{CGGR}. The discrete dual process is easier than the original
continuous process and can be used to study detailed time dependent properties of the wealth distribution.\\

In section 5, we introduce a class of diffusion processes naturally associated to the wealth distribution model, in the same way in which the KMP model is related to the Brownian energy process by instantaneous thermalization, see \cite{GKRV}.\\

The two-agent models that we study here are to be considered as the building block of more complex models, because models with $N$ agents are always based
on transactions between two agents, or put otherwise, in the language of kinetic theory, one restricts to binary collisions.
Moreover, the results that we obtain for two agents models remains partly true and contain relevant information for more complex $N$ agent models: the absence of product
invariant measures for non-zero propensity, 
characterization of the redistribution laws (of beta-type) for which there are product measures
in the zero propensity case, and the existence of a discrete dual process.

Moreover, the diffusion
of expected wealth as a random walk 
where
the propensity only enters in the time scale is actually shown here in section 6 in the $N$ agent model with general symmetric transactions.
This shows in particular that these models satisfy ``Fourier law'', and have a ``finite positive conductivity''.
In the two agent case we can moreover completely characterize the set of stationary distributions, also in the presence
of non-zero propensity. The latter becomes more complex in models with more agents, and because complex
dependencies will appear, it is not expected that for general models we can obtain analytic results. In forthcoming work we will
however
analyze the covariances between different agents in more complex $N$ agents models and show how they are related
to the Green's function of the underlying random walk.
\color{black}

An important subject in the area of wealth distribution models is the emergence of power laws \cite{eli3,eli2,eli}, and in particular
of the Pareto law for the stationary wealth distribution marginals. In our models, Pareto laws only emerge in a particular choice
of the ``redistribution law''. It is expected that in more complex models with $N$ agents and {\em annealed random
propensities}, power laws can emerge due to the fact that with non-zero probability many agents have propensity close
to one. Such a result is actually obtained in \cite{DM} for a particular model with deterministic redistribution law.

Our paper is organized as follows: Section 2 introduces the general models for energy redistribution and wealth distribution and the relative notation; Section 3 is devoted to the study of stationary measures in these models; in Section 4 duality is introduced and discussed; Section 5 deals with the diffusions associated to the two sites energy redistribution models and finally Section 6 analyzes the case of a simple economy with $N$ agents, by studying how the expected wealth of a single agent diffuses.

\section{Models with two agents}
\subsection{Notation}
Let $E=[0,\infty)^2$ be the state space of the models with two agents. An element $(x,y)\in E$ represents
the wealth of agent 1 resp.\ 2, or in the energy redistribution model the energy at vertex 1, resp.\ vertex 2.

For $\epsi\in [0,1]$ we define the map $T_\epsi:E\to E$ as
\be\label{teps}
T_\epsi (x,y) = (\epsi(x+y), (1-\epsi)(x+y)),
\ee
and, for $\lambda, \epsi\in [0,1]$, the map $T_\epsi^\lambda:E\to E$ as
\begin{eqnarray} \label{tepsla}
T_\epsi^\lambda (x,y)&=& \left(\lambda x+ (1-\lambda)\epsi(x+y), \lambda y+ (1-\lambda)(1-\epsi) (x+y)\right)\\ \nonumber &=&\lambda(x,y) +(1-\lambda)T_\epsi(x,y).
\end{eqnarray}
Notice that the map $T^\lambda_\epsi$
conserves the total energy $s=x+y$. Naturally $T_\epsi= T^0_\epsi$, whereas $T^1_\epsi$ is the identity.
\bd
A redistribution measure $\nu= \nu(s,d\epsi)=\nu(s,\epsi) d\epsi$ is a weakly continuous map
from $[0,\infty)$ into the space of absolutely continuous probability measures with full support on $[0,1]$.
\ed
The interpretation of $\nu(s,d\epsi)$ is the distribution of the fraction $x/(x+y)$
of the total wealth going to the first agent after one redistribution
when the total wealth is equal to $s=x+y$ before (and hence also after) the redistribution.

\subsection{The energy redistribution model}
For a given redistribution measure $\nu= \nu(s,d\epsi)$, the energy distribution model
based on $\nu$
is a Markov process in which, for an initial state $(x,y)\in E$, after an exponential waiting time (with mean one),
the values $x,y$ are replaced by $T_\epsi(x,y)$, where the redistribution parameter $\epsi$ is distributed according to $\nu(s, d\epsi)$.\\
The generator of the energy redistribution model is defined on bounded continuous functions $f:E\times E\to\R$ as
\be\label{engen}
\loc f(x,y)= P f(x,y)- f(x,y),
\ee
where the one-step transition redistribution operator $P$ is given by
\be\label{redesp}
P f(x,y)= \int_0^1 f(T_\epsi (x,y))\nu(x+y, \epsi) d\epsi.
\ee

\br
Important special cases are obtained when the redistribution measure $\nu(x+y, d\epsi)$ does not depend on $s=x+y$.
\begin{itemize}
\item[a)]
When $\nu(x+y,d\epsi)= d\epsi$ is uniform, the corresponding redistribution model is the well-known KMP model, introduced in \cite{KMP} as a microscopic model of heat conduction
satisfying Fourier law.
\item[b)]
If $\nu(x+y,d\epsi)=\nu(d\epsi)$ is the Beta distribution with support  $[0,1]$ and  parameters $(m,m)$, i.e.
$\nu(x+y,d\epsi)= \epsi^{m-1} (1-\epsi)^{m-1} \frac1{B(m,m)} d\epsi$
(where $B$ denotes the Beta function), then the redistribution model coincides with the so-called thermalized Brownian energy process studied in \cite{CGGR}.
\end{itemize}
In these cases, the stationary measures are of product form. For more details about the existence of product measures we refer to section 3 below.
\er
\subsection{The wealth distribution model}
In a wealth distribution model the pair $(x,y)$ represents the wealth of
two agents, and redistributions are the results of random transactions between these two agents. These agents represent the whole economy, whose total wealth is preserved. The interest for these type of models is 
mainly in trying to understand wealth transactions at the micro level. In that sense, considering just two agents is only the starting point for analyzing a more realistic economy made of $N$ agents that interact among them. Two is in fact the minimum number of agents we need to model a transaction. As we will partially see later on, knowing what happens at the level of single pair is fundamental to make inference on the whole $N$ agents economy.

Compared with the energy redistribution model, the wealth distribution model has an additional parameter $\lambda\in (0,1)$ interpreted as the
propensity to save of the agents. This parameter represents the fraction of wealth that the agents do not exchange; upon each redistribution only the wealth $(1-\lambda)(x+y)$ is transacted among the two agents.\\
For a given redistribution measure $\nu (x+y, \epsi)d\epsi$ and $\lambda\in (0,1)$ the
wealth distribution model with parameters $(\nu,\lambda)$ has the generator
\be\label{wealthgen}
\loc^\lambda f(x,y)= P^\lambda f(x,y)- f(x,y),
\ee
acting on bounded continuous functions, and
where the redistribution kernel is now equal to
\be\label{redeswp}
P^\lambda f(x,y)= \int_0^1 f(T^{\lambda}_{\epsi}(x,y))\nu(x+y, \epsi)d\epsi.
\ee
\br
In the literature (see e.g.  \cite{CCCC}), different models are available with $\epsi$ fixed, $\lambda$ random, or both $\epsi$ and $\lambda$ random.\\
In this paper we restrict our attention to the case in which $\epsi$ is random and $\lambda$ is fixed.\\
Another more realistic variant is an agent-dependent propensity to save $\lambda_{1,2}$, in which the map $T^{\lambda}_\epsi$ becomes
\begin{equation}
\label{l1l2}
T^{\lambda_1,\lambda_2}_\epsi (x,y) = (\lambda_1 x + \epsi((1-\lambda_1) x+ (1-\lambda_2) y), \lambda_2 y + \epsi((1-\lambda_1) x+ (1-\lambda_2) y).
\end{equation}
\er

\subsection{Questions}
The basic questions we focus on in this paper are the following:
\ben
\item What are the stationary distributions?
\item When do we have product stationary distributions?
\item Do these models have discrete dual processes, and if so what are the corresponding duality functions? (See Section \ref{se4} for more details.)
\een
Answering these questions is meaningful for the development of more realistic wealth distribution models, in which $N$ agents, possibly belonging to different classes, interact.\\
Regarding question 1 it is important to distinguish between {\em canonical} and {\em grand canonical invariant} measures. \\
Since for all these models $s=x+y$ is a conserved quantity, upon fixing $s$, the process
will converge to a unique invariant measure concentrating on the set $\{(x,y): x+y=s\}$. This invariant measure is called the canonical measure, and especially in the two agents case, these measures are quite easy to find out. 
All invariant measures are mixtures of these canonical measures. The canonical measures are of course never product, and do
not have full support, because
the sum $x+y$ is fixed. Mixtures of canonical measures, i.e., choosing a particular distribution for $s=x+y$ can however
be product. Mixtures of canonical measures, where we allow $s=x+y$ to fluctuate are called grand canonical measures.

In generalizing from the $2$ agents case to the case of $N$ agents, it is important to find out whether these canonical measures correspond to grand canonical measures with full support (i.e. where $s=x+y$ is not fixed), and in which cases the latter can be product measures.
We will see that for non-zero propensity $\lambda>0$, there do not exist mixtures of canonical measures
that are product, which implies that there do not exist product stationary distributions in that case.

\section{Stationary measures}
\subsection{Product measures for the energy redistribution model}
The following theorem gives a complete characterization of the cases in which the energy redistribution model generates product measures.
\bt\label{prodthm}
Let us consider an energy redistribution model, then:
\ben
\item
The product measure $\mu(x)\mu(y) dxdy$ is invariant if and only if
\be\label{brom}
\nu(s, a) = \frac{\mu(as) \mu((1-a)s)}{\int_0^1\mu(\alpha  s) \mu((1-\alpha)s) d\alpha}.
\ee
This measure is then also reversible. More generally $\mu(x,y) dx dy$ is invariant if and only if
\be\label{bram}
\nu(s, a) = \frac{\mu(as, (1-a)s)}{\int_0^1\mu(\alpha  s,(1-\alpha)s) d\alpha}.
\ee
\item
If $\nu(s,\epsi)=\nu(\epsi)d\epsi$ does not depend on $s$, then the only invariant product measure $\mu$ is of Gamma type, i.e.
\[
\mu(x)= \lambda^{2b}\frac{x^{2b-1}}{\Gamma( 2b-1)} e^{-\lambda x}.
\]
The corresponding Beta-distribution for $\nu$ is equal to
\[
\nu(\epsi)= \frac{1}{B(b,b)}\epsi^{b-1} (1-\epsi)^{b-1}.
\]
\een
\et

\bpr
We start from point (1). If $\mu(x)\mu(y) dxdy$ is invariant then for all $f:E^2\to\R$ bounded and continuous we have
\[
\int Pf (x,y) \mu(x)\mu(y) dx dy= \int f(x,y) \mu(x)\mu(y) dx dy,
\]
which, using the change of variables $a= x/(x+y), s=x+y$ leads to the following equalities:
\begin{eqnarray*}
&&\int_0^\infty\int_0^\infty f(x,y) \mu(x)\mu(y) dx dy
\\
&=& \int_0^\infty \int_0^1 f(as, (1-a)s) \mu(as)\mu((1-a)s) sds da
\\
&=&
\int Pf(x,y) \mu(x)\mu(y) dx dy
\nonumber\\
&=&
\int_0^\infty\int_0^\infty\int_0^1 f(\epsi(x+y), (1-\epsi)(x+y)) \mu(x)\mu(y)\nu(x+y,\epsi) dx dy d\epsi
\\
&=&
\int_0^\infty\int_0^1\int_0^1 f(\epsi s, (1-\epsi)s) \mu(as)\mu((1-a)s)\nu(s,\epsi) s ds da d\epsi
\\
&=&
\int_0^\infty\int_0^1 f(a s, (1-a)s) \left(\int_0^1\mu(\epsi s)\mu((1-\epsi)s) d\epsi \right) \nu(s,a)s ds da.
\end{eqnarray*}
This proves \eqref{brom}, since we infer that for every bounded continuous $f$
\begin{eqnarray*}
&&\int_0^\infty \int_0^1 f(as, (1-a)s) \mu(as)\mu((1-a)s) sds da\\
&=& \int_0^\infty\int_0^1 f(a s, (1-a)s) \left(\int_0^1\mu(\epsi s)\mu((1-\epsi)s) d\epsi \right) \nu(s,a)s ds da.
\end{eqnarray*}
If \eqref{brom} is fulfilled, then with a similar computation, one obtains
\[
\int Pf(x,y) g(x,y) \mu(x)\mu(y) dxdy= \int Pg(x,y) f(x,y) \mu(x)\mu(y) dxdy,
\]
which implies that $\mu$ is also reversible.
When $\mu$ is not product, the computation is identical. This finishes the proof of item 1.

For item (2), if $\nu(s, a)$ does not depend on $s$, then the equation for
$\mu$ reads
\be\label{booo}
\nu(a) = \frac{\mu(as)\mu((1-a)s)}{\int_0^1 \mu(\alpha s) \mu((1-\alpha)s) d\alpha}.
\ee
By symmetry of the rhs of \eqref{booo} under $a\leftrightarrow (1-a)$, $\nu$ has to be of the form $K(a) K(1-a)$. Setting $\psi(s)=\int_0^1 \mu(\alpha s) \mu((1-\alpha)s) d\alpha$, we find
\be\label{kwak}
\frac{\mu(as)\mu((1-a)s)}{K(a)K(1-a)}=\psi(s).
\ee
Putting $J(a,s)= \mu(as)/K(a)$ and differentiating \eqref{kwak} w.r.t. $a$, we get
\[
\frac{\partial_a J(a,s)}{J(a,s)}= \frac{\partial_a J(1-a,s)}{J(1-a,s)},
\]
where $\partial_a$ denotes partial derivative w.r.t.\ $a$. This leads to
\[
s\frac{\mu'(as)}{\mu(as)} - \frac{K'(a)}{K(a)}= s\frac{\mu'((1-a)s)}{\mu((1-a)s)} + \frac{K'(1-a)}{K(1-a)}.
\]
Setting $a=0$, we now have
\[
\frac{\mu'(s)}{\mu(s)} = C_1 + \frac{C_2}{s},
\]
which leads to $\mu(s)= s^c e^{-c' s}$. Since we require $\mu(s)$ to be normalized on $[0,\infty)$,
we must have $c'>0$, which gives us a Gamma distribution.
\epr
\\
In the following theorem we show the general structure of the invariant measures for the energy redistribution model.
\bt\label{invthm}
Every distribution of the type
$(\epsi S, (1-\epsi)S)$, where $S$ is a non-negative random variable, and
$\epsi$ is distributed according to the redistribution measure $\nu(s, d\epsi)$, is invariant for the energy redistribution
model with generator $\loc$. These measures exhaust the set of invariant measures.
\et
\bpr
Let $\mu$ be the joint distribution $(\epsi S, (1-\epsi)S)$ and denote by $\gamma$ the distribution of $S$.
Then we have by assumption that $\mu(dx dy) =\gamma (ds) \nu(s, d\epsi)$. Thus
\beq
\int P f(x,y) \mu(dx dy)&=&
\int\left(\int_0^1 f(\epsi(x+y),(1-\epsi)(x+y)) \nu(s,d\epsi)\right)\mu(dxdy)
\nonumber\\
&=&
\int_0^1 \int_0^\infty\int_0^1 f(\epsi s, (1-\epsi)s) \nu(s,d\epsi) \gamma(ds)  \nu(s,d\epsi')
\nonumber\\
&=&\int_0^1 \int_0^\infty f(\epsi s, (1-\epsi)s) \nu(s,d\epsi) \gamma(ds)
\nonumber\\
&=&\int f(x,y) \mu(dxdy)
\eeq
Therefore, $\mu$ is invariant.

For fixed $s=x+y$, the invariant distribution $\mu_s$ concentrating on the pairs $(x,y)$ with $x+y=s$ is unique and
ergodic, and hence given
by the distribution of $(\epsi s, (1-\epsi)s)$, where $\epsi$ is distributed according to the redistribution measure
(corresponding to $\gamma= \delta_s$.
These are the only ergodic invariant measures. Since the set of invariant measures forms
a simplex, every invariant measure is of the form $\int \gamma (ds) \mu_s$.
\epr
\br
In the special cases in which $\nu$ is both independent of $s$ and a Beta distribution in
$\epsi$, the only choice for a distribution $\gamma$ on the $s$-variable, so that the couple $(\epsi s, (1-\epsi)s)$ is i.i.d.,
is the Gamma distribution. All this shows consistency between theorems \eqref{invthm} and \eqref{prodthm}.\\
If we drop the requirement of full support, and choose the redistribution measure as
\[
\nu(s, \epsi)=\epsi^{-\alpha-1}(1-\epsi)^{-\alpha-1}, \qquad \frac1s\leq \epsi \leq \frac{s-1}{s},
\]
then the corresponding product measure is the product of Pareto I distributions, i.e.
\[
\mu(x,y) =\mu(x)\mu(y),
\]
with
\[
\mu(x) =\alpha x^{-\alpha-1}, \qquad x\geq 1.
\]
\er

\subsection{Stationary measures for the wealth distribution model}\label{se4}

Let us define a new variable $r= x/s$, representing the proportion of wealth owned by agent $x$. As usual $s=x+y$ is the total wealth. Under a Markov process with generator
\eqref{wealthgen}, we have that $s_t=s_0=s$ and $r_t$ is Markovian with $s$-dependent generator
\be\label{les}
L^s f(r)= \int (f(\lambda r+ (1-\lambda)\epsi)-f(r))\nu(s, d\epsi).
\ee
If $\nu$ does not depend on $s$ this process has an $s$-independent generator
\be\label{less}
Lf(r)= \int(f(\lambda r+ (1-\lambda)\epsi)-f(r))\nu( d\epsi).
\ee
The unique stationary distribution of this process is the
unique distribution that is stationary under the recursion
\be\label{recuo}
r_{n+1}= \lambda r_n + (1-\lambda) \epsi_n,
\ee
where $\epsi_n$ is an i.i.d. sequence with distribution $\nu$.\\
This distribution can be easily obtained by iterating the recursions \eqref{recuo}, and it is thus given by the distribution of
\be\label{epsla}
\epsi^\lambda_\infty= \sum_{n=0}^\infty \lambda (1-\lambda)^n \epsi_n,
\ee
where the different $\epsi_n$ are i.i.d. with distribution $\nu$.\\
In the following theorem we characterize the stationary measure of the
process with generator \eqref{wealthgen} and $s$-independent $\nu$. For a redistribution
measure $\nu$ we denote by $\nu^\lambda_\infty$ the distribution of the corresponding
random variable $\epsi^\lambda_\infty$.
\bt
Consider the wealth distribution model with generator \eqref{wealthgen} and $\nu(s,\epsi)= \nu(\epsi)$. Then all invariant measures are distributions of the form
\[
(\epsi^\lambda_\infty S, (1-\epsi^\lambda_\infty)S),
\]
where $S$ is a non-negative random variable, and where $\epsi^\lambda_\infty$ is given in \eqref{epsla}. For $\lambda>0$ there are no product invariant measures.
\et
\bpr
Let $\mu$ be an invariant measure and call $\gamma(ds)$ and $\lambda(s,dr)$
the corresponding distributions of $s=x+y$ and $r=x/s$. Then it follows that $\lambda(s,dr)$ must be the unique invariant measure of
the process with generator \eqref{les}, which is $\nu^\lambda_\infty$.\\
Such a measure can only be product if
\[
\gamma(x+y) \nu^\lambda_\infty (x/(x+y))= \psi(x)\psi(y).
\]
Then, by the same reasoning as in the proof of point (2) of theorem \ref{prodthm}, we conclude that
$\gamma(x)$ must be of the form $x^c e^{-c'x}$, while $\nu (r)=C r^b (1-r)^b$, which is a contradiction, because $\nu^\lambda_\infty(a)$ is clearly not of this form.
\epr
\br
If $\nu$ does depend on $s$, then for each $s$ we have a corresponding
$\nu^\lambda_\infty(s, \cdot)$, and the stationary measures are given by
$(\epsi_\infty S, (1-\epsi_\infty)S)$, where {\em given $S=s$},
$\epsi_\infty$ is distributed according to $\nu^\lambda_\infty(s)$.
\er

\br
If the propensity is agent dependent (cf. \eqref{l1l2}), then the stationary measures can still be found explicitly.
In this case the analogue of the recursion \eqref{recuo} reads
\be\label{recuop}
r_{n+1}= r_n ( \lambda_1 + (\lambda_2-\lambda_1) \epsi_n) + (1-\lambda_2) \epsi_n
\ee
which leads to the limiting distribution given by
\be\label{twoproprinfty}
r_\infty= \sum_{k=0}^\infty (1-\lambda_2) \epsi_k \prod_{i=0}^{k-1}( \lambda_1 + (\lambda_2-\lambda_1)\epsilon_i)
\ee
where the empty product is defined to be equal to $1$. The corresponding stationary measures are the distributions of $(r_\infty S,(1-r_\infty) S)$, where $S$ is a non-negative random variable.
\er

\section{Duality}

\subsection{Duality for the energy redistribution model: $s$-independent case}

Duality is an important tool in the study of interacting particle systems, since it allows to connect two Markov processes via a duality function. In other words it makes possible to build a connection between processes characterized by a continuous state space, such as the energy and wealth distribution models (or processes of diffusion type), and discrete processes, without any loss of information. Given that discrete processes are generally easier to handle, both analytically and for the
purpose of simulations, duality represents a very powerful instruments for the study of many complex models.
\bd
Let $\{\xi_t :t\geq 0\}, \{\eta_t :t\geq 0\}$ be two Markov processes with state spaces $\Omega$ and  $\Omega'$ respectively. Let $D:\Omega\times\Omega'$ be a measurable function.
Set $\widehat{E}_\xi$ to denote expectation in the $\xi$ process starting from $\xi\in\Omega$, and the same for $\E_\eta$ with $\eta\in \Omega'$. \\
We say that the process $\{\xi_t:t\geq 0\}$ is a dual process of the process $\{\eta_t:t\geq 0\}$
with duality function $D$ if, for all $\xi\in\Omega$ and for all $\eta\in \Omega'$,
we have
\be\label{dualityrel}
\widehat{E}_\xi( D(\xi_t, \eta))= \E_\eta (D(\xi, \eta_t)).
\ee
We denote this property by using the following notation
\[
\{\xi_t:t\geq 0\}\longrightarrow^D \{\eta_t:t\geq 0\}.
\]
\ed
The most interesting use of duality arises when the dual process is simpler than the original one, for example discrete instead of continuous, or with finitely many particles
instead of infinitely many, and when the duality functions are ``rich enough'' to preserve ``full information''. \\
A classical example is the choice $D(n,x)= x^n$ as a candidate duality function between a process with a continuous state space variable $x\in \R$ and a discrete process $n\in\N$ usually of birth and death type.\\
Duality was discovered and used for the KMP model in \cite{KMP}, and for more general energy redistribution models in \cite{CGGR}.\\
We will now show that the energy redistribution process with generator
as in equation \eqref{engen}, and with
$\nu(s,d\epsi)$ independent of $s$, {\em always admits a discrete dual process}, with duality functions of the form
\be\label{dualform}
D(n,m; x,y) = \frac{x^n y^m}{c_{nm}}
\ee
However, {\em only in the case in which $\nu$ is a Beta distribution}, these duality functions will {\em factorize}, that is to say
$c_{nm}= d_n d_m$. In that case, the dual process is known as the thermalized
$SIP$ (symmetric inclusion process), introduced and studied in \cite{CGGR}.\\
Let us now start by acting with the generator \eqref{engen} on $f_{nm} = x^n y^m$. Denoting
$\nu_{nm}=\int_0^1 \epsi^n (1-\epsi)^m \nu(\epsi) d\epsi$, we obtain
\begin{eqnarray*}
\loc f_{nm} &=&
\left(\sum_{k=0}^{n+m} {{n+m}\choose k} \nu_{nm} f_{k,n+m-k}\right) - f_{n,m}
\\
&=&
\sum_{r,s: r+s=n+m} a_{nm;rs} f_{rs}.
\end{eqnarray*}
The coefficients $a_{nm;rs}$ are non-negative when $(n,m)\not= (r,s)$. We can then
apply a general result from \cite{EBG}:
if $\mu$ is an invariant measure and
\begin{equation} \label{cnm}
c_{nm} = \int x^n y^m \mu(dx\ dy),
\end{equation}
then by setting
\be\label{dualq}
q_{nm, rs}= a_{nm, rs} \frac{c_{rs}}{c_{nm}},
\ee
we have that the $q$'s define the transition rates of a continuous-time Markov chain on $\N\times\N$, and we have the duality between
the generators
\be\label{gendual}
\loc D(n,m; x,y) = \sum_{rs} q_{nm, rs} D(r,s; x,y),
\ee
with $D(n,m;x,y)$ given by \eqref{dualform}.
Equation \eqref{gendual} implies that, in terms of the associated processes,
\be\label{dualpr}
\E_{x,y} (D(n,m, X_t,Y_t) ) = \hat{\E}_{n,m} (D(N_t, M_t; x,y)),
\ee
that is to say
\[
\{ (N_t, M_t):t\geq 0\} \longrightarrow^D \{ (X_t, Y_t):t\geq 0\}.
\]
The dual Markov chain can be interpreted as a discrete redistribution model, where the total number of ``dual particles'' $(n+m)$ is conserved and ``dual particles'' are redistributed
according to the transition rates $q_{nm, rs}$.\\
Coming back to the simplest case where $\nu(d\epsi)= d\epsi$, we find $\nu_{nm}= \frac{n!m!}{(n+m+1)!}$. For the invariant product measure we choose
$\mu(dx dy)= e^{-x} e^{-y} dx dy$, $c_{nm}= n! m!$, and we find
\[
q_{nm, rs}= \delta_{n+m, r+s}  \frac{1}{n+m+1},
\]
i.e. the dual particles are uniformly redistributed. This is exactly the dual KMP model,
or thermalized SIP with $m=2$ \cite{KMP}, \cite{CGGR}. \\
Similarly, when $\nu$ is a Beta distribution,
we find a hypergeometric redistribution rule for the dual rates, corresponding to a thermalized $SIP$.\\
These are the only cases in which the coefficients $c_{nm}$ factorize, because factorization occurs if and only if the invariant measure $\mu$ in \eqref{cnm} is a product measure. Therefore, by Theorem \ref{prodthm}, the coefficients $c_{nm}$ factorize if and only if $\nu$ is
a Beta distribution.\\
We can then summarize our findings in the following theorem.
\bt
The energy distribution model with $s$-independent redistribution measure admits a
discrete dual Markov jump process with duality functions
\[
D(n,m;x,y)= \frac{x^n y^m}{c_{nm}}.
\]
The coefficients $c_{nm}$ factorize if and only if the
redistribution measure is Beta distributed, that is when $\nu(\epsi)= \frac1{B(b,b)}\epsi^{b-1} (1-\epsi)^{b-1}$. In that case
\[
c_{nm} =\frac{\Gamma(b)}{\Gamma( b+n)}\frac{\Gamma(b)}{\Gamma(b+m)}
\]
\et
Notice that the coefficients $c_{nm}$ in the duality function depend on the choice of the invariant measure
$\mu$. So choosing one reference invariant measure
$\mu_0$ for the choice of the coefficients, i.e.,
$c_{nm} = \int x^n y^m \mu_0(dx dy)$,  we have that if $\mu$ is a (possibly different) invariant measure for the process, then
defining
\[
\int D(n,m; x,y) \mu(dx dy)= \widehat{\mu} (n,m)
\]
we have, as a consequence of invariance of $\mu$
that $\widehat{\mu}$ is harmonic for the dual process, i.e.,
\[
\widehat{\E}_{(n,m)} (\widehat{\mu} (n_t,m_t))= \widehat{\mu} (n,m)
\]
As a consequence, since the dual chain is irreducible once $n+m$ is fixed, the only harmonic functions of the dual process are functions of
the conserved quantity $n+m$. Therefore, we obtain that
\be\label{baka}
\frac{\int x^n y^m \mu(dx dy)}{c_{nm}}= \phi (n+m)
\ee
i.e., the joint moments under $\mu$ are the same as under $\mu_0$ up to
a $n+m$ dependent factor. Conversely, if \eqref{baka} holds, then $\mu$ is also
an invariant measure.
\subsection{Duality for the wealth distribution model: $s$-independent case}

Let us suppose once again that $\nu(s,d\epsi)=\nu(d\epsi)$. It is convenient to study the duality
of the wealth distribution model
via the process $r_t$ with generator \eqref{less}.
Let us consider the function $f_n (r)= r^n$. Acting on the generator, we find
\be\label{predual}
L f_n = \sum_{k=0}^n a(n,k) f_k,
\ee
where for $k\not= n$
\be\label{acoe}
a(k,n)=  {n\choose k} \lambda^k (1-\lambda)^{n-k} \int \epsi^{n-k}\nu(d\epsi),
\ee
and $a(n,n)= \lambda^n -1 <0$. Therefore, the conditions of
the general theorem in \cite{EBG} are satisfied, and
we have a dual process $(n_t: t\geq 0)$ jumping
from
$k$ to $n$, $k<n$ at rate
\be\label{dualrate}
q(k,n)= \frac{\alpha_k}{\alpha_n} {n\choose k} \lambda^k (1-\lambda)^{n-k} \int \epsi^{n-k}\nu(d\epsi),
\ee
with
\[
\alpha_r = \E((\epsi^\lambda_\infty)^r),
\]
and with duality functions
\be\label{dualfu}
D(n,r) = \frac{r^n}{\alpha_n}.
\ee
If we want to consider duality between the $\{X_t, Y_t\}_{t \geq 0}$ process and a discrete process, it
is convenient to introduce the following function
\[
K(n_1,n_2, r)= r^{n_1}(1-r)^{n_2}.
\]
Acting with the generator \eqref{less}, we then find a dual process
$\{n_1(t), n_2 (t)\}$ going from
$\{n_1, n_2\}$ to $\{k_1, k_2\}$ at rate
\be\label{ra}
Q(k_1, k_2, n_1, n_2)= A(k_1, k_2, n_1, n_2)\frac{\alpha(k_1,k_2)}{\alpha(n_1,n_2)},
\ee
where
\be\label{rb}
\alpha(i,j)= \int r^i (1-r)^j \nu^\lambda_\infty (dr),
\ee
and where
\be\label{rc}
A(k_1, k_2, n_1, n_2) = {n_1\choose k_1} {n_2\choose k_2}\lambda^{k_1+k_2} (1-\lambda)^{n_1+n_2-k_1-k_2} \int \epsi^{n_1-k_1}(1- \epsi)^{n_2-k_2}\nu(d\epsi).
\ee
The duality function is then given by
\be\label{rd}
D(n_1, n_2, r) = \frac{r^{n_1}(1-r)^{n_2}}{\alpha(n_1,n_2)}
\ee
If $x+y=1$, then $r=x$ and $1-r= y$, and we can rewrite the duality function as
\be\label{buro}
D(n_1, n_2, r) = \frac{x^{n_1}y^{n_2}}{\alpha(n_1,n_2)}.
\ee
As a consequence, we obtain the following duality result.

\bt
Between the processes $\{(X_t, Y_t), t\geq 0\}$, with $X_0+Y_0=1$ with generator \eqref{wealthgen} and
$s$ independent $\nu$, and the discrete jump process $\{(n_1(t), n_2(t)):t\geq 0\}$,
with rates given by \eqref{ra}, \eqref{rb}, \eqref{rc}, we have
the duality relation
\[
\{(n_1(t), n_2(t)):t\geq 0\}\longrightarrow^D \{(X_t, Y_t), t\geq 0\}
\]
with $D$ given in \eqref{buro}.
\et

\section{Diffusions on two sites associated to energy redistribution models}

The energy redistribution models can be considered as ``thermalization limits'' of slower
processes where the energy is redistributed continuously, as in the case of energy diffusion processes. We investigate these diffusions in the present context, where they can be considered as models of wealth distribution on more microscopic time scales.\\
The concept of thermalization limit
has been introduced in \cite{GKRV} to build a connection between the Brownian energy process,
a process where energy diffuses, and the KMP model, which is an energy redistribution model (with uniform redistribution measure).\\
The idea of thermalization is that at exponential waiting times a pair of vertices is selected, and the diffusion process associated
to that pair is ``run for infinite time'', so that the energies associated to that pair of vertices
are redistributed according to the stationary distribution of the diffusion process.\\
 In this section we go the other way round: the energy or wealth distribution model is given, and
we look for the corresponding diffusion process, for which this model is the thermalization limit.\\
Since for two sites $s=x+y$ is conserved, for the diffusion process we can start with a generator of the type
\be\label{difgen}
L f(x,y) = xy(\partial_x-\partial_y)^2 + a(x,y) (\partial_x-\partial_y),
\ee
where $\partial_x$,$\partial_y$ denote partial derivatives with respect to $x$,$y$.\\
This is the generator of a diffusion process $\{(X_t, Y_t): t\geq 0\}$ on $[0,\infty)^2$ conserving the sum $s=x+y$.
Using $s=x+y$ and $r= x/(x+y)$, we obtain that under
the diffusion $(X_t, Y_t)$ we have $s_t=s_0$ is conserved and $r_t$
evolves according to the diffusion process on $[0,1]$ (with absorbing boundaries)
with generator
\be\label{epsidif}
L_s^r= r(1-r) \partial_r^2 + \frac{a(sr, (1-s)r)}{s} \partial_r
\ee
This diffusion has a stationary distribution $\psi(r) dr$
with density
\be\label{statmeas}
\psi (r) = C(1-r)^{-1} r^{-1} \exp\left(\int\frac{a(r s, (1-r)s)}{s} dr \right)
\ee
with $C$ a normalizing constant.\\
Let us first consider the simple linear case $a(x,y) = -\alpha(x-y)$, with $\alpha>0$. Then the diffusion for $r$
becomes $s$-independent and its generator
reads
\[
L_s^r= r(1-r) \partial_r^2 + \alpha(2r-1)\partial_r
\]
and the corresponding stationary measure
\[
\psi(r) = C(1-r)^{\alpha-1} r^{\alpha-1}
\]
i.e., the $Beta(\alpha,\alpha)$ distribution.

More generally, if we want $\nu(s,r)$ to be the stationary distribution of the generator
$L_s^r$ then by solving \eqref{statmeas} for $a$ we arrive at the following relation between $\nu$ and $a$:
\be\label{anu}
a(r s, (1-r)s) = sr (1-r) \frac{d}{dr} (\log ((1-r)r \nu(s,r))).
\ee
Note that once more inserting here $\nu(s,r)= C (1-r)^{\alpha-1} r^{\alpha-1}$ consistently
gives $a(x,y)= -\alpha(x-y)$.

The connection between the diffusion model and the energy redistribution model is then
as follows. Since the diffusion has a unique stationary measure $\nu(s, dr)$ for $x+y=s$ fixed
\[
\lim_{t\to\infty}\E_{x,y}f(X_t,Y_t) = \int f(r s, (1-r)s)\nu(s, dr),
\]
i.e. the generator of the energy redistribution model can be written as
\[
\loc f(x,y) = \lim_{\beta\to\infty} (e^{\beta L} f(x,y)- f(x,y)),
\]
i.e. it is the thermalization of the diffusion generator.

\section{Random walk of the expected wealth in the wealth distribution model with $N$ agents}

Although for the wealth distribution model with $N$ agents one cannot in general easily characterize
the stationary distribution (which will certainly not be product because of the
results of the two site model), a general first-order duality result can be nevertheless obtained: the expectation of $x_i(t)$ can be described via a single random walk.\\
To define the wealth distribution model with $N$ agents, we consider
a countable set of agents denoted by $S$ (with cardinality $N$), and a symmetric
random walk transition probability $p(i,j)$ indexed by $i,j\in S$, i.e. $p(i,j)=p(j,i)\geq 0$ and $\sum_{j\in S} p(i,j)=1$, for all $i\in S$.\\
We denote by ($\caX_t:t\geq 0$) the continuous random walk, jumping
at rate one according to the transition probabilities $p(i,j)$, and by $\caE_i$ the expectation in this process starting from $i\in S$.\\
The state space is then $E= [0,\infty)^S$. For $x= \{x_i, i\in S\}\in E $ and $i,j\in S$, we denote by $T_{\lambda, \epsi}^{ij}(x)$ the configuration obtained from $x$ by replacing the couple $(x_i,x_j)$ with $ T^{\lambda}_{\epsi}(x_i,x_j)$, while leaving all other coordinates as they are.\\
Given a ($s$-independent) redistribution measure $\nu$ on $[0,1]$, the wealth distribution model
has the following generator, defined on bounded continuous functions depending on a finite number of coordinates,
\be\label{genwealthgen}
\loc^\lambda f(x) = \sum_{ij\in S} 2p(i,j) \int(f(T_{\lambda, \epsi}^{ij}x)- f(x))\nu(x+y, d\epsi).
\ee
Let us denote by $\E_x$ the expectation in this process starting from $x\in E$. We further assume that $\int \epsi \nu(x+y, d\epsi)=1/2$. Hence we have the following general result.

\bt
In the wealth distribution model with generator \eqref{genwealthgen}, we have
\be\label{lindual}
\E_{x}( x_i(t)) = \caE_i (x_{\caX((1-\lambda)t)})= \sum_{j} p_{(1-\lambda )t}(i,j) x_j,
\ee
where $p_t(i,j)$ is the transition probability of the random walk $\caX_t$.\\
As a consequence, if $\mu$ is an invariant measure and
\[
\rho_i = \int x_i \mu(dx),
\]
then $\rho$ is harmonic for $p(i,j)$, i.e.
\[
\sum_{j\in S} p(i,j) \rho_j = \rho_i.
\]
\et
\bpr
For $f_i(x)= x_i$, using the symmetry of $p(i,j)$, i.e.
\begin{eqnarray*}
\left(\loc^\lambda f_i \right)(x)&=&
\sum_j p(i,j) \left(\lambda x_i +(1-\lambda)\int\epsi \nu(x+y, d\epsi) (x_i+x_j)-x_i\right)
\\
&+&
\sum_j p(j,i) \left(\lambda x_i + (1-\lambda)\int(1-\epsi) \nu(x+y, d\epsi) (x_i+x_j)-x_i\right)
\\
&=& \sum_{j} p(i,j)(1-\lambda) (f_j(x)-f_i(x))
\\
&=&
(1-\lambda) \left(A f_{\text{\bf{\Large{.}}}}(x) \right)(i),
\end{eqnarray*}
with $A$ the generator of the process $\caX_t$, and $f_{\text{\bf{\Large{.}}}}(x)$ the function $i\mapsto f_i(x)$.\\
This shows that
$\E_x( f_i(x(t)))=:\psi(t,i, x)$ satisfies
\[
\frac{d}{dt} \psi(t,i,x)= (1-\lambda)(A\psi(t, \cdot, x))(i).
\]
This equation has the unique solution $\psi(t,i,x)=\caE_i (\psi(0, \caX_{t(1-\lambda)}, x))$.\\
This proves \eqref{lindual}. To prove
the second statement, it is sufficient to integrate \eqref{lindual} over a stationary measure $\mu$.
\epr

\section{Conclusion}
In this paper we analyze a two-agent model of energy and wealth distribution, serving as
building block for more complex and realistic models. These models are parametrized by a redistribution law
$\nu$
and a saving propensity $\lambda$.
We show that non-zero saving propensity necessarily leads to dependencies in the stationary distribution.
We show that only for beta-distributed redistribution laws $\nu$ 
we can obtain product stationary distributions. Moreover the models always admit a discrete
dual process and also we introduce a class of naturally associated diffusion processes.

For a general $N$ agent model with symmetric redistribution we demonstrate that the expected wealth moves
as a random walk and the propensity only enters as a time scale. These results are
first steps in the understanding of more complex $N$ agents models, including emergence
of power laws and long-range dependencies.

\end{document}